\documentclass[12pt]{article}

\usepackage{graphicx}
\usepackage{setspace} 

\input epsf.sty
\input psfig.sty

\def\comment#1{}
\newcommand{\remove}[1]{}
\newcommand{\ignore}[1]{}

\begin{document}
\bibliographystyle{plain}
%\doublespacing
\title{Condorcet Winner Probabilities - A Statistical Perspective}
\author
{ M.\ S.\ Krishnamoorthy\\
Department of Computer Science,\\Rensselaer Polytechnic Institute, Troy, NY 12180\\
and\\
M.\ Raghavachari\\ 
Decision Sciences and Engineering Systems,\\ Rensselaer Polytechnic Institute, Troy, NY 12180
}
\date{}

\maketitle
\begin{newpage}
\begin{abstract}
A Condorcet voting scheme chooses a winning candidate as one who defeats all
others in pairwise majority rule. We provide a review which includes the 
rigorous mathematical
treatment for calculating the limiting probability of a Condorcet winner
for any number of candidates and value of $n$ odd or even and with arbitrary rank order probabilities,
when the voters are independent. We provide a compact and complete Table for
the limiting probability of a Condorcet winner with three candidates and
arbitrary rank order probabilities. We present a simple proof of a result
of May to show the limiting probability of a Condorcet
winner tends  to zero as the number of candidates tends to infinity.
We show for the first time that the limiting probability of a
Condorcet winner for any given number of candidates $m$ is monotone decreasing in
$m$ for the equally likely case. This, in turn, settles the conjectures of 
Kelly and Buckley and Westen for the case $n \rightarrow \infty$.
We prove the validity of Gillett's conjecture on the minimum value
of the probability of a Condorcet winner for $m=3$ and any $n$. We generalize
this result for any $m$ and $n$ and obtain the minimum solution and the
minimum probability of a Condorcet winner.

\end{abstract}
\begin{newpage}

\section{Introduction}
Different voting procedures (see for example \cite{mue03}  Page 147) for electing a leader or deciding
a candidate or a judgment have been proposed in the literature. In all
voting procedures, we assume that the voter rank orders the candidate or
decision from highest to lowest. One of the most used schemes is a majority 
rule,
in which a candidate or decision is ranked first by more than half the number of voters. An example of the majority rule is the U.S. Supreme Court decisions.
Another selection strategy is the plurality rule in which a candidate is
ranked first by the largest number of voters. An example of the plurality rule
is the U.S. presidential election. In the Condorcet scheme, a candidate is
chosen who defeats all others in pairwise majority rule. An example of the
Condorcet scheme is the selection of job candidates or proposals to be funded
in academic departments and by funding agencies or for determining winners in tournaments.

While the Condorcet scheme is a transitive in individual voter's choices, 
 there are many instances in which
a majority winner may not exist \cite{nie68}\cite{may71}. For example, let there be 3
voters ($V_1,V_2,V_3$) and 3 candidates ($C_1,C_2,C_3$). $V_1$ rank 
orders $C_1C_2C_3$, $V_2$ rank orders $C_2C_3C_1$ and finally $V_3$ 
rank orders $C_3C_1C_2$. In this voting example, $C_1$ is preferred 
over $C_2$ two out of three times. Similarly $C_2$ is preferred 
over $C_3$ and $C_3$ is preferred over $C_1$ two
out of three times. Hence there is no one with majority preferences 
over all pairs. So according to the Condorcet scheme, there is no winner in this
case.

In past research work, it has been assumed that a probability distribution is
given to all the rank orders. It is also assumed in most studies that all voters vote independently.
Previous researchers have calculated the probability of a Condorcet winner using
analytical and simulation models \cite{jon95} by making simplifying assumptions on the
probability distributions. Asymptotic expressions for the limiting probability of
a Condorcet winner (i.e., when the number of voters $\rightarrow$ $\infty$) 
 have been obtained \cite{gil52} \cite{geh98} \cite{gar68} \cite{nie68} \cite{may71}.  A survey of previous research in this area may be
found in \cite{geh83}.

The focus of the present paper is on the limiting probability of a Condorcet
winner when the number of voters, $n$, tend to infinity for a given number $m$ of
candidates and arbitrary probability distribution on the $m!$ rankings. 
We present a review and some new results in this part of the study. While
researchers in this area, notably for example, Niemi and Weisberg \cite{nie68},
Garman and Kamien \cite{gar68} have noted the relationship of the limit with
certain multivariate normal probabilities, there doesn't exist a statistically
rigorous proof of this result in order to justify the evaluation of the limit
for various cases. Further it is often assumed that the number of candidates 
is odd, see Niemi and Weisberg \cite{nie68}. We close this gap with this
paper by providing a rigorous statistical proof for any number of candidates, 
the number of voters both odd and even, 
and arbitrary probability distributions, usually referred to as culture
probabilities. In section 2, we show further that for any number of candidates
 and with arbitrary probability distributions on the preference ranking, the
limiting probability depends on the calculation of positive orthant
probabilities of appropriate multinormal distributions.
We also provide a 
compact and complete
Table covering all possible cases for computing the limiting probability of a Condorcet winner when the
number of candidates is 3 and with arbitrary probability distribution among the
rank orders. We show that the limiting probability as $n\rightarrow \infty$ for
any number of candidates $m$ is monotone decreasing in $m$ which partially
validates Kelly \cite{kel74} and Buckley and Westen \cite{bucwes79} conjectures. 
We also present a simple proof of May's result \cite{may71} 
that the limiting
probability of a Condorcet winner for equally likely case tends to zero
as the number of candidates gets larger. We also prove Gillett's \cite{gil80}
conjecture on the minimum probability of a Condorcet number for 3 candidates
and any number of voters. We extend the result to any $m$ and $n$ and obtain 
the minimum value and solution. In Section 5, we treat the case of
$m=4$ candidates with arbitrary probabilities $p_i$ and show that it is
possible to have exact expressions for the limiting Condorcet winner
probability for all possible scenarios.

\section{Condorcet Winner: Statistical Formulation and Results for General Case}
In this section, we present a statistical derivation of the limiting 
probability of a Condorcet winner
when the voters are independent and the number of voters is large. Let
\begin{center}
\begin{tabular}{ccl}
$m$ &=& Number of Candidates\\
$n$&=& Number of Voters \\\
\end{tabular}
\end{center}
We have $K~=~m!$ preference rankings of $m$ candidates. A voter will choose one
of these rankings. We assume that the voters act independently. Let $p_i$ be 
the probability that a voter prefers the rank order $i$, for $i~=~1,2, \cdots , K$. 
Further, we know that $p_i ~ \ge ~ 0$ and $\displaystyle\sum_{i=1}^{K} p_{i}= 1$.  
Let $N_i$ be the number of voters voting for the $i$th preference ranking, for $i=1,2, \cdots , K$. We therefore have $\displaystyle\sum_{i=1}^{K} N_{i} = n$. Now, it is 
well known that $(N_1,N_2, \cdots ,N_K)$ has a multinomial distribution.
\begin{equation}
P(N_i=n_i, i=1, \cdots, K)= \frac{n!}{n_{1}!n_{2}! \cdots n_{K}!} \prod_{i=1}^{K} p_{i} ^{n_i}
\label{eqn1}
\end{equation} 
where the $n_i$'s are nonnegative integers with $\displaystyle\sum_{i=1}^{K}n_i=n$.

It is known that (see, for example,\cite{cra48} Page 318), for $i=1, \cdots , K$, $j=1, \cdots , K$ and $ i \ne j$.
\begin{eqnarray*}
E[N_i] &=& n p_i,\\
Var[N_i] &=& n p_i (1-p_i),\\
Covariance(N_i,N_j)&=& -n p_i p_j, \\
Correlation(N_i,N_j)&=& -\sqrt \frac {p_i p_j} {(1-p_i)(1-p_j)}
\end{eqnarray*}

Denote the $m$ candidates by $C_1,C_2, \cdots , C_m$. Write $C_iPC_j$, $i \ne j$ if $C_i$ is preferred to $C_j$ by a voter.
Define $C_iMC_j$, $i \ne j$ if $C_i$ beats $C_j$ by majority rule. i.e., the
total number of voters preferring $C_i$ to $C_j$ is more than the total number
of voters preferring $C_j$ to $C_i$. This fact can be mathematically expressed as
\begin{equation}
a_{1,(i,j)} N_1+a_{2,(i,j)} N_2 + \cdots + a_{K,(i,j)} N_K \ge 1
\label{eqn2}
\end{equation}
where the $a_{l,(i,j)}$'s are $\pm$ 1. $a_{l,(i,j)}~=~1$, if in the $l$th preference
ranking voted by the $N_l$ voters, $C_iPC_j$. Similarly $a_{l,(i,j)}~=~-1$ if
$C_jPC_i$ is the preference ranking voted by $N_l$ voters. 
Clearly $\frac{K}{2}$ of the $a_{l,(i,j)}$'s are 1 and $\frac{K}{2}$ of 
the $a_{l,(i,j)}$'s are equal to -1. 

Equation (\ref{eqn2}) implies that the difference between the number of voters
preferring $C_i$ to $C_j$ and the number of voters preferring $C_j$ 
to $C_i$ is at least 1. We define, for example, $C_i$ to be the Condorcet winner if and only
if $C_iMC_j$ for each $j~=~1,2, \cdots , m$ and $j \ne i$. For example, the probability that $C_1$ is the
Condorcet winner is the probability of the joint event
\begin{eqnarray}
a_{1,(1,2)} N_1 + a_{2,(1,2)} N_2 + \cdots + a_{K,(1,2)} N_K & \ge & 1 \nonumber \\
a_{1,(1,3)} N_1 + a_{2,(1,3)} N_2 + \cdots + a_{K,(1,3)} N_K & \ge & 1 \nonumber \\
\cdots & \ge & 1 \nonumber \\
\cdots & \ge & 1 \nonumber \\
a_{1,(1,m)} N_1 + a_{2,(1,m)} N_2 + \cdots + a_{K,(1,m)} N_K & \ge & 1 \label{eqn3}
\label{eqnsw}
\end{eqnarray}
Condition (\ref{eqnsw}) lead to what is defined as 'strong' winner. Replacing
all the 1's by 0's lead to the definition of a 'weak' winner. For small values $n$ and $m$ , this probability can be computed using the 
multinomial distribution (\ref{eqn1}). We are interested in finding \\
\(  \lim_{n \rightarrow \infty} \) $P(\hbox{One of the $C_i$ is a Condorcet winner} )$
for a given value of $m$.

Several researchers have attempted to find this limit or an approximation to it.
Gilbaud \cite{gil52} gave an expression of this limit for $m=3$, when all the $p_i$'s are
equal (Impartial Culture or IC model), though he did not explain the method. 
Garmen and Kamien \cite{gar68} derived the expression for the limiting Condorcet
winner probability for $m=3$ and $m=4$ for the Impartial Culture model without
giving details.
Jones et al \cite{jon95},  
Bell \cite{bel70} and Niemi et al \cite{nie68} obtained approximations to it by extensive simulations or quadrature methods attributed to Ruben \cite{rub54} 
for
the Impartial Culture (IC) model i.e., when all the preference rankings 
are equal.

Niemi and Weisberg \cite{nie68} and Garman and Kamien \cite{gar68} were the first to relate the limit
or an approximation to it with the calculation of certain multivariate normal
probabilities. 
Niemi et al \cite{nie68} assume $n$ to be an odd number and Garman et al 
\cite{gar68} indicate the relation in a footnote. 
However,
they do not present a rigorous statistical proof that the limit is exactly 
a multivariate normal
probability for any general culture probabilities. They state that the multivariate normality is achieved by the
fact that the multinomial distribution tends to the multivariate normal
distribution. The multivariate normality however, is the consequence of the
fact that a set of linear combinations of $N_i$ values, 
defined by equation (\ref{eqn2}), tends to the multivariate normal distribution
by the application of multivariate central limit theorem, as we proceed to show
here. In this section we show a rigorous statistical derivation of the
limit which covers the cases when n is odd or even.

Most researchers in this area have focused on the odd $n$ (number of voters)
case. This is because ties do not occur while determining the majority rule; i.e.,
for odd number of voters $n$, the number of voters preferring $C_i$ to $C_j$
will never equal the voters preferring $C_j$ to $C_i$. However, if $n$ is
even there is a positive probability of having a tie. We show, in what 
follows, that the probability of a tie for $n$ even tends to zero as $n~ \rightarrow~\infty$.
Let $p_{ij}=\displaystyle\sum_r { p_r}$ the with summation extending over the preference
rankings in which $C_i$ is preferred over $C_j$ and let $T$ denote the total
number of voters preferring $C_i$ to $C_j$.
Then, it is well known from Multinomial distribution theory that
$T$ is a Binomial variable with $n$ and $p_{ij}$ as parameters. The probability of a tie in
the majority determination between $C_i$ and $C_j$ is
\begin{eqnarray}
P(T=\frac{n}{2})&=& {n \choose {\frac {n}{2}}} {p_{ij}}^{\frac{n}{2}} {(1-p_{ij})}^{\frac{n}{2}} \nonumber \\
&=& \frac{n!}{{(\frac {n}{2})!} {(\frac {n}{2})!}}{[p_{ij}(1-p_{ij})]}^{\frac{n}{2}} \label{eqn31}
\end{eqnarray}

We know by Stirling's approximation \cite{cra48} Page 130, that
$$
n!~\cong \sqrt{{2\pi}}e^{-n} n^{n+\frac{1}{2}}
$$
The symbol $\cong$ means $\frac{n!}{\sqrt{2\pi}e^{-n}n^{n+\frac{1}{2}}} \rightarrow 1$ as $n \rightarrow \infty$. Applying this to equation (\ref{eqn31}) and
noting that $p_{ij}(1-p_{ij}) \le \frac{1}{4}$ for all $0\le p_{ij}\le 1$, we
can verify that 
\begin{eqnarray*}
P(T=\frac{n}{2}) &\cong& \frac{1}{\sqrt{{2\pi}}} \frac{2^{n+1}} {\sqrt{n}} {[p_{ij}(1-p_{ij})]}^{\frac{n}{2}} \\
&\le& \frac{1}{\sqrt{2n\pi}}\frac{2^{n+1}}{4^{\frac{n}{2}}}\\
&=& \sqrt{\frac{2}{\pi}} \frac{1}{\sqrt{n}} \rightarrow 0~ \hbox{as}~ n\rightarrow \infty 
\end{eqnarray*}

For a given voter $v$, for $i=1,\cdots,m$, j$=1,\cdots,m$, $i~\ne~j$, define
\begin{eqnarray}
X(i,j,v) &=&1 \hbox{~~if~~ $C_iPC_j$} \nonumber \\
&=&-1~~ \hbox{~~if~~ $C_jPC_i$} \label{eqn4} 
\end{eqnarray}
For the voter profile $(N_1,N_2,\cdots, N_K)$, $C_iMC_j$ if 
$\displaystyle\sum_{v=1}^{n} X(i,j,v) \ge 1 $ . The above representation has been used in Bell (\cite{bel70}) and
Gehrlein (\cite{geh98}). Note that
\begin{eqnarray}
\displaystyle\sum_{v=1}^{n} X(i,j,v)& =& a_{1,(i,j)} N_1 + a_{2,(i,j)} N_2 + \cdots + a_{K,(i,j)} N_K  \nonumber \\
\hbox{We have,} \nonumber \\
\lambda_{ij} &=& E(X(i,j,v)) \nonumber \\
~~&=& a_{1,(i,j)} p_1 + a_{2,(i,j)} p_2 + \cdots + a_{K,(i,j)} p_K \label{eqn49} \\
Var(X(i,j,v))&=& 1~-~\lambda_{ij}^2 \label{eqn5}
\end{eqnarray}

We also need to obtain the correlation coefficient between 
$X(i,j,v)$ and $X(i,l,v)$, $i \ne j$, $i \ne l$ and $j \ne l$. To express correlation, we first
define
\begin{eqnarray}
a_{r,(i,j),l} &=& 1 ~~\hbox{if in the preference ranking corresponding to $N_r,r=1,\cdots,K$ is} \nonumber\\
&~&~~\hbox{$C_iPC_j$ and $C_iPC_l$}\nonumber\\
&~&~~\hbox{or $C_jPC_i$ and $C_lPC_i$}\nonumber \\
&=& -1 ~~\hbox{Otherwise} \label{eqn7}
\end{eqnarray}
Covariance is defined by the following equation.
\begin{eqnarray*}
Covariance(X(i,j,v),X(i,l,v))&=&E[X(i,j,v)X(i,l,v)]-\lambda_{ij}\lambda_{il}
\end{eqnarray*}
The correlation between $X(i,j,v)$ and $X(i,l,v)$ is given by
\begin{eqnarray*}
R_{jl}^{(i)}&=& \frac {E(X(i,j,v).X(i,l,v))~-~\lambda_{ij}\lambda_{il}} {\sqrt{(1-\lambda_{ij}^2)(1-\lambda_{il}^2)}} ~~\hbox{$j=1,\cdots,m$,}\\ 
&&~~\hbox{$l=1,\cdots,m$ and $j \ne l$, $j \ne i$, $l \ne i$.} \\
&=&\frac{ \displaystyle\sum_{r=1}^{K} {a_{r,(i,j),l}p_r -\lambda_{ij}\lambda_{il}}} {\sqrt{(1-\lambda_{ij}^2)(1-\lambda_{il}^2)}}
\end{eqnarray*}

For a given $i=1,2,\cdots,m$, let $R_i$ denote the $(m-1) \times (m-1)$ 
correlation matrix $( R_{jl}^{(i)})$. Note that $R_{jj}^{(i)}=1, j=1,2,\cdots,m ; j \ne i$.

Let
\begin{eqnarray}
y(i,j,v)&=& \frac{(X(i,j,v)-\lambda_{ij})}{\sqrt{1-\lambda_{ij}^2}}, j=1,2,\cdots,m, j \ne i \nonumber \\
\hbox{and~~~} Z_n(i,j)&=& \displaystyle\sum_{v=1}^n {\frac{y(i,j,v)}{\sqrt{n}}}\nonumber\\
&=&\displaystyle\sum_{v=1}^n{\frac{(X(i,j,v)-\lambda_{ij})}{\sqrt{n}\sqrt{1-\lambda_{ij}^2}}}, j=1,\cdots,m, j \ne i \label{eqn6}
\end{eqnarray}
Note that the $(m-1)$ random variables $y(i,j,v),j=1,\cdots,m, j \ne i$ have a
joint distribution with zero mean vector, unit variances and Correlation matrix $R_i$.

From (\ref{eqn6}), the probability that $C_iMC_j$ is given by
\begin{equation}
P\left(\displaystyle\sum_{v=1}^{n} {X(i,j,v) \ge 1}\right)\\
\Longrightarrow P\left(Z_n(i,j)\ge \frac{(\frac{1}{\sqrt{n}}-\sqrt{n}\lambda_{ij})}{\sqrt{1-\lambda_{ij}^2}}\right).
\label{eqn71}
\end{equation}

The probability that $C_i$ is the Condorcet winner is therefore given by the joint
probability
\begin{eqnarray}
P\left[Z_n(i,j)\ge \frac{(\frac{1}{\sqrt{n}}-\sqrt{n}\lambda_{ij})}{\sqrt{1-\lambda_{ij}^2}}, \hbox {$j=1,\cdots,m$, ~$j \ne i$ }\right]. \label{eqn8}
\end{eqnarray}

By the multivariate central limit theorem (See Cram\`er \cite{cra70}, Page 112) applied to independent 
summands $Z_n(i,j)$,$j=1,\cdots,m$, $j \ne i$, the joint distribution 
of 
$$
(Z_n(i,1),\cdots,Z_n(i,i-1),Z_n(i,i+1),\cdots,Z_n(i,m))
$$ 
tends as $n \rightarrow \infty $ to the $(m-1)$ dimensional multivariate normal distribution with 
zero mean vector, unit variances and correlation matrix $R_i~=~(R_{jl}^{(i)})$. 
Let us write
$$
\delta_{ij} = \left\{ \begin{array}{rl}
                      -\infty &\mbox{if $\lambda_{ij} >0$}\\
                      \infty &\mbox{if $\lambda_{ij} <0$} \\
                      0 &\mbox{if $\lambda_{ij} =0$} \\
		      \end{array}
              \right.
$$
Denote by $L(h_1,h_2,\cdots,h_{m-1};R)$ the multivariate normal probability
$$ 
P(Z_1^* \ge h_1, Z_2^* \ge h_2, \cdots , Z_{m-1}^* \ge h_{m-1};R)
$$
where $(Z_1^*,Z_2^*,\cdots,Z_{m-1}^*)$ has the multivariate normal distribution
with zero mean vector, unit variances and Correlation matrix $R$. Then we have
$$
\lim_{n \rightarrow \infty}P\left[Z_{ij} > \frac{(\frac{1}{\sqrt{n}}-\sqrt{n}\lambda_{ij})}{\sqrt{1-\lambda_{ij}^2}}, \hbox {$j=1,\cdots,m$, ~$i \ne j$ }\right]
$$
\begin{equation}
=L(\delta_{i1},\delta_{i2},\cdots,\delta_{i,i-1},\delta_{i,{i+1}},\cdots,\delta_{im};R_i).
\label{eqnl}
\end{equation}
The limiting probability of a Condorcet winner is then
\begin{equation}
P(\infty,m)~=~\displaystyle\sum_{i=1}^m{L(\delta_{i1},\delta_{i2},\cdots,\delta_{i,i-1},\delta_{i,{i+1}},\cdots,\delta_{im};R_i)}.
\label{eqnll}
\end{equation}

Equation (\ref{eqnll}) shows that the limiting Condorcet winner probability is obtained
as the sum of $m$ $L$-functions. The arguments of the $L$ functions are either
$+\infty$ or $-\infty$ or 0. If at least one of the arguments in a $L$
function is $+\infty$, the function value is 0 by the definition of the
$L$-function. We can also drop all the arguments with value $-\infty$ and calculate the reduced $L$-function from the marginal multivariate normal
distribution. It is therefore interesting to note that even for arbitrary
probability distribution $p_i$'s, at worst we have to evaluate positive
orthant probabilities of a multivariate normal distribution. We shall give
examples for the case $m=3$ and $m=4$ in later sections and provide a 
scheme to obtain the exact limiting probability expressions for a Condorcet winner
for any choice of $p_i$'s.

Equation (\ref{eqnll}) can also be specialized to the Dual Culture(DC) Model.
A preference rank order is the dual of another ranking if each can be 
obtained by reversing the order of the other. For example, if $m=4$,
rank order $C_1C_4C_3C_2$ is dual to $C_2C_3C_4C_1$. The DC model assigns
$p_i$ such that $p_i=p_j$ if $i$ is dual to $j$. In this case it is clear from
the definition of $\lambda_{ij}$ that all the $\lambda_{ij}  \equiv 0$. This, 
in turn, implies that $\delta_{ij} \equiv 0$. Equation (\ref{eqnll}) shows 
then that the limiting probability of a Condorcet winner under DC models
$$
\displaystyle\sum_{i=1}^m L(0,0,0,\cdots,0;R_i)
$$

Note that in each $L$ function, there are $(m-1)$ zeroes. Gehrlein \cite{geh97} gives the exact expressions for $m=3$ and $m=4$.

\section{Condorcet Winner: Some Special Cases}
First we assume that all preference rankings for a voter are equally likely.
This implies that $p_i~\equiv~1/K$, for $i=1,\cdots,K$. For this case,
symmetry conditions imply that $\lambda_{ij}~\equiv~0$ for $i=1,\cdots,m$ ,
$j=1,\cdots,m$ and $i\ne j$ $\Rightarrow \delta_{ij}~ \equiv~0$. For this case,
it can be verified from \cite{geh76} that $R_i~\equiv ~R^*$ where $R^*$ is 
the $(m-1) \times (m-1)$ equicorrelated matrix with equal correlation 
value = 1/3. The limiting probability of a Condorcet winner from (\ref{eqnll})
 is given by
\begin{eqnarray}
P(\infty,m)&=&mL(0,0,\cdots,0;R^*) \label{eqn9}
\end{eqnarray}

This expression for $m=3$ was obtained by Gilbaud (\cite{gil52}) as 
$\frac{3}{4}+\frac{3}{2\pi} \arcsin(\frac{1}{3})$. 
Garman and Kamien \cite{gar68} gave the expression for $m=4$ without
giving details. For $m \ge 5$, several researchers have obtained approximate values for expression (\ref{eqn9}) based on simulation
studies or by quadrature methods. See for example Gehrlein \cite{geh76}, Jones et al \cite{jon95} and Niemi et al \cite{nie68}.
The expression (\ref{eqn9}) shows that the limiting probability for $m$ candidates can be
obtained by an evaluation of positive orthant probabilities for a $(m-1)$ dimensional
multivariate normal distribution with zero mean vector, unit variances and equicorrelated correlation matrix with all correlations equal to 1/3. Several papers
in the statistical literature deal with such evaluations. See, for example, David and Mallows \cite{dav61}, Sondhi\cite{son61}, Placket\cite{pla54}, Ruben\cite{rub54}, 
Johnson and Kotz \cite{jon72} and Bacon \cite{bac63}. Bacon \cite{bac63} derives explicit expressions for the
positive orthant probabilities for a few small values of $m$ for a multivariate
normal distribution with equicorrelated correlation matrix with all correlations
equal to $\rho$. He also derives a recursive relation to calculate the probabilities for successive values of $m$. Except for a few small values of $m$, such expressions involve the evaluation of
multiple integrals without known closed form expressions. Denote the 
expression (\ref{eqn9}) by $P(m)$. Using the expressions from Bacon 
\cite{bac63}, Gehrlein \cite{geh83} obtains 
the values of $P(m)$ for some small values of m as

\begin{tabular}{ll}
$m=3$,&$P(3)=\frac{3}{4} + \frac{3}{2\pi} \arcsin(\frac{1}{3})$\\
$m=4$,&$P(4)=\frac{1}{2}[1+\frac{6}{\pi} \arcsin(\frac{1}{3})]$\\
$m=5$,&$P(5)=\frac{5}{16}[1+\frac{12}{\pi} \arcsin(\frac{1}{3}) + \frac{24}{{\pi}^2} \int_0^{\frac{1}{3}} {\arcsin (\frac{\lambda}{1+2\lambda})\frac{d\lambda}{\sqrt{(1-\lambda^2)}}} ]$\\
$m=6$,&$P(6)=\frac{3}{16}[1+\frac{20}{\pi} \arcsin(\frac{1}{3}) + \frac{120}{{\pi}^2} \int_0^{\frac{1}{3}} {\arcsin (\frac{\lambda}{1+2\lambda})\frac{d\lambda}{\sqrt{(1-\lambda^2)}}} ]$\\
$m=7$,&$P(7)=\frac{7}{64} [1+\frac{30}{\pi} \arcsin(\frac{1}{3}) +\frac{360}{\pi^2}\int_0^{\frac{1}{3}} {\arcsin (\frac{\lambda}{1+2\lambda})\frac{d\lambda}{\sqrt{(1-\lambda^2)}}}  $\\
&$+\frac{720}{\pi^3} \int_0^{\frac{1}{3}}\int_0^{\frac{\mu}{1+2\mu}} {\arcsin (\frac{\lambda}{1+2\lambda})\frac{d\lambda}{\sqrt{(1-\lambda^2)}} \frac{d\mu}{\sqrt{(1-\mu^2)}} } ]$\\
\end{tabular}

Note however, that some of these expressions involve integrals.

Denote $L(0,0,\cdots,0,\rho)$ by $L_{m-1}(\rho)$ when there are $m-1$ zeros in
the argument of $L$ and $\rho$ is the common correlation. From the results of  
Sampford (See Moran \cite{mor52}) we have
$$
L_{m-1}(\rho)=\frac{1}{\sqrt{\pi}} \int_{-\infty}^{\infty} { e^{-t^2}{(1-\Phi(at))}^{m-1} dt}
$$
where $\Phi(.)$ is the cumulative distribution function of the standard
normal distribution and $a=\frac{2\rho}{1-\rho}$. 
For the IC model, $\rho=\frac{1}{3}$ and hence $a=1$. Thus

\begin{equation}
L_{m-1}(0,\cdots,0;\frac{1}{3})=\frac{1}{\sqrt{\pi}} \int_{-\infty}^{\infty} {e^{-t^2} {(1-\Phi(t))}^{m-1} dt}
\label{eqnaa}
\end{equation}

It is also known that $L_{m-1}(\rho)$ for
$m$ being even can be obtained recursively from $L_i(\rho)$ for $i \le m-2$:
\begin{equation}
L_{m-1}(\rho)=\frac{1}{2}\left[\frac{-(m-1)}{2}+1+{{m-1} \choose 2}L_2(\rho)-{{m-1} \choose 3}L_3(\rho)+\cdots+{{m-1} \choose {m-2}} L_{m-2}(\rho)\right]
\label{eqnab}
\end{equation}
See Johnson and Kotz \cite{jon72} Pages 48-52 for more details. We obtained the
expression for $m=6$ using the recursion formula.

For $m > 7$, Bacon \cite{bac63} provides an approximation believed to be reasonably 
close to the true value. Figure \ref{fig:asymp} shows the plot of $P(m)$
 for increasing values of $m$ based on these approximations. 
May \cite{may71} has shown that $P(m) \rightarrow 0$, as 
$m \rightarrow \infty$ at the rate of $\frac{1}{m}$.
His proof is based on saddle point approximation of an integral
and it appears to be complicated and difficult to understand. We
present a direct proof of this result.

From (\ref{eqn9}) and (\ref{eqnaa}), we have
\begin{eqnarray*}
P(m)&=& \frac{m}{\sqrt{\pi}} \left[\int_{-\infty}^{\infty} {e^{-t^2} {(1-\Phi(t))}^{m-1} dt}\right] \\
&=&\sqrt{2}m \left[\int_{-\infty}^{\infty} {\phi(t){(1-\Phi(t))}^{m-1} e^{\frac{-t^2}{2}} dt}\right] 
\end{eqnarray*}
where 
$$
\phi(t)=(\frac{1}{\sqrt{2\pi }}) e^{\frac{-t^2}{2}} = \frac{d}{dt} \Phi(t)
$$
\begin{eqnarray*}
P(m)&=&\sqrt{2}m\left[\left\{\frac{-{(1-\Phi(t))}^m e^{\frac{-t^2}{2}}}{m}\right\}_{-\infty}^{\infty}-\frac{1}{m} \int_{-\infty}^{\infty} {{[1-\Phi(t)]}^m t e^{\frac{-t^2}{2}} dt}\right]\\
&=&-\sqrt{2} \int_{-\infty}^{\infty} {{[1-\Phi(t)]}^m e^{\frac{-t^2}{4}} t e^{\frac{-t^2}{4}} dt}
\end{eqnarray*}
Applying Cauchy-Schwarz's inequality 
\begin{equation}
\left|{\int {f(t)g(t) dt}}\right| \le {\left(\int {f^2(t) dt}\right)}^{\frac{1}{2}} {\left(\int {g^2(t) dt}\right)}^{\frac{1}{2}}\nonumber \\ 
\end{equation}
to the above equation, we get
\begin{eqnarray*}
P(m) &\le& \sqrt{2}{\left(\int_{-\infty}^{\infty} {{(1-\Phi(t))}^{2m}e^{\frac{-t^2}{2}} dt}\right)}^{\frac{1}{2}} {\left(\int_{-\infty}^{\infty} {t^2 e^{\frac{-t^2}{2}} dt}\right)}^{\frac{1}{2}}
\end{eqnarray*}

Note that
\begin{eqnarray*}
\int_{-\infty}^{\infty} {{(1-\Phi(t))}^{2m} e^{\frac{-t^2}{2}} dt} &=& \sqrt{2\pi} \int_{-\infty}^{\infty}{{(1-\Phi(t))}^{2m}\phi(t) dt}\\
&=&\frac{\sqrt{2\pi}}{2m+1}\\ \\
\int_{-\infty}^{\infty} {t^2 e^{\frac{-t^2}{2}} dt} &=& \sqrt{2\pi} \int_{-\infty}^{\infty}{t^2\phi(t) dt}\\
 &=&\sqrt{2\pi}
\end{eqnarray*}

Thus 
$$
P(m) \le 2\pi\sqrt{2} \frac{1}{{(2m+1)}^{\frac{1}{2}}} \rightarrow~0 \hbox{~as~} m \rightarrow \infty
$$

We proved that the limiting probability $P(m)$ of a Condorcet winner tends to 
zero as $m \rightarrow \infty$. We show in addition that$P(m)$ is 
monotonically decreasing in $m$. This settles Kelly's \cite{kel74} and Buckley and 
Westen's \cite{bucwes79} conjecture for large values of $n$. This monotonicity property, we
believe, is new.

{\bf Proposition:} $P(m+1) < P(m)$ \\

{\bf Proof:} From Equations (\ref{eqn9}) and (\ref{eqnaa}), we have to show 
that

\begin{equation}
\frac{(m+1)}{\sqrt{\pi}}\int_{-\infty}^{\infty} {e^{-t^2} {(1-\Phi(t))}^{m} dt}
<  \frac{m}{\sqrt{\pi}}\int_{-\infty}^{\infty} {e^{-t^2} {(1-\Phi(t))}^{m-1} dt}
\label{eqnA} 
\end{equation}
\begin{equation}
\frac{(m+1)}{\sqrt{\pi}}\int_{-\infty}^{\infty} {e^{-t^2} {(1-\Phi(t))}^{m} dt}
 =  \sqrt{2} (m+1) \int_{-\infty}^{\infty} {e^{\frac{-t^2}{2}} {(1-\Phi(t))}^{m} \phi(t) dt}
\label{eqnB}
\end{equation}
where $\phi(t)$ = pdf of standard normal distribution. The right side of
equation (\ref{eqnB}) equals after integrating by parts

$$
-\sqrt{2} \int_{-\infty}^{\infty}  {t e^{\frac{-t^2}{2}} {(1-\Phi(t))}^{m+1} dt}
$$

Similarly the right side of equation (\ref{eqnA}) equals
$$
-\sqrt{2} \int_{-\infty}^{\infty}  {t e^{\frac{-t^2}{2}} {(1-\Phi(t))}^{m} dt}
$$

Thus it is enough to show 
$$
\int_{-\infty}^{\infty}  {t e^{\frac{-t^2}{2}} \left( {(1-\Phi(t))}^{m+1} -{(1-\Phi(t))}^{m} \right)  dt} >0
$$
Split the region of integration into $\int_{-\infty}^0 + \int_0^{\infty} $,put
$t=-y$ in the first integral and after simplification and noting $1-\Phi(-y)=\Phi(y)$,
it reduces to show that
\begin{equation}
 \int_{0}^{\infty}  { \left[ \left( {(1-\Phi(y))}^{m+1} -{(1-\Phi(y))}^{m} \right) - \left(\Phi^{m+1}(y)-\Phi^m(y) \right)\right] y \phi(y)  dy} >0
\label{eqnC}
\end{equation}

We have written $\Phi^{m+1}(y)$ for ${\Phi(y)}^{m+1}$ and $\Phi^{m}(y)$ for ${\Phi(y)}^{m}$.
We will show that
\begin{equation}
\left( \Phi^{m}(y) -{(1-\Phi(y))}^{m} \right) - \left(\Phi^{m+1}(y) -{(1-\Phi(y))}^{m+1}\right) >0 \hbox{~~for $y>0$}
\label{eqnD}
\end{equation}

Write $\Phi(y)=z$ and note that the left hand side of equation (\ref{eqnD}) equals
\begin{equation}
(2z-1)\left[z^{m-1}+z^{m-2}(1-z)+\cdots+(1-z)^{m-1}\right]-(2z-1)\left[z^m+z^{m-1}(1-z)+\cdots+(1-z)^m\right]
\label{eqnE}
\end{equation}

Observe that $2z-1>0$ for $y>0$.

Note that $z^{m-1}-z^{m-1}(1-z) = z^m$. Further,
$$
z^{m-2}(1-z)-z^{m-2}(1-z)^2=z^{m-1}(1-z)
$$
and so on.

Hence the expression in the square brackets of equation (\ref{eqnC})
$$
= z^{m-1}(1-z)+z^{m-2}(1-z)^2+\cdots + z{(1-z)}^{m-1} \\
>0
$$
Thus the integrand in equation (\ref{eqnC}) is $>$0 and hence the left side of 
equation (\ref{eqnC}) is $>0$.

\begin{figure}[tb]
\begin{center}
\begin{tabular}{c}
\setlength{\epsfxsize}{4.5in} {\epsfbox{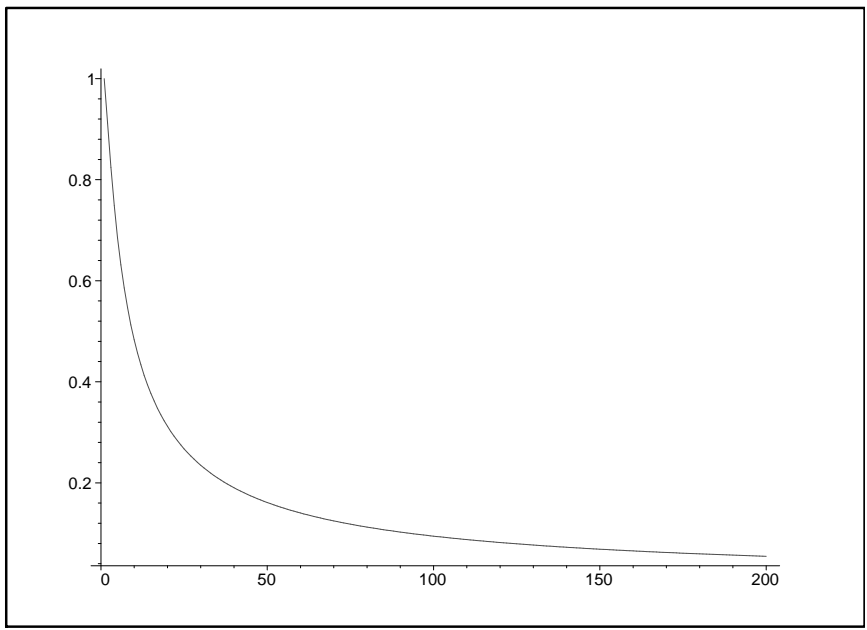}}
\end{tabular}
\end{center}
\caption{The plot of limiting Probability of a Condorcet winner vs the number of Candidates}
\label{fig:asymp}
\end{figure}

\section {Condorcet Winner: Special Case $m=3$ candidates}

Niemi et al \cite{nie68} calculated the limiting probability of a Condorcet winner 
values for $n=1..49$, utilizing the table of Ruben \cite{rub54}.
Ruben used quadrature methods and recursive formulas to obtain these values. 
May \cite{may71} and Garman et al \cite{gar68} discussed the case of 3 candidates to obtain the limiting
probability of a Condorcet winner. Gehrlein \cite{geh97} discusses several
possibilities that may arise in obtaining the value for $m=3$. We give a compact and complete
Table for all possible $p_i$'s utilizing the results of Section 2 and (\ref{eqnll}) of this paper.
Suppose the 6 rankings are \\
\begin{tabular}{|l|l|l|}
\hline \hline 
Rankings&Number of Voters&Probabilities\\ \hline \hline
$C_1C_2C_3$&$N_1$&$p_1$\\ \hline 
$C_1C_3C_2$&$N_2$&$p_2$\\ \hline 
$C_2C_1C_3$&$N_3$&$p_3$\\ \hline 
$C_2C_3C_1$&$N_4$&$p_4$\\ \hline 
$C_3C_1C_2$&$N_5$&$p_5$\\ \hline
$C_3C_2C_1$&$N_6$&$p_6$\\ \hline 
\end{tabular}
\\
\\
\\
\\
\noindent

It can be verified that
\begin{eqnarray*}
\lambda_{12}&=&p_1+p_2+p_5-p_3-p_4-p_6 \\
\lambda_{13}&=&p_1+p_2+p_3-p_4-p_5-p_6 \\
\lambda_{21}&=&p_3+p_4+p_6-p_1-p_2-p_5 \\
\lambda_{23}&=&p_1+p_3+p_4-p_2-p_5-p_6 \\
\lambda_{31}&=&p_4+p_5+p_6-p_1-p_2-p_3 \\
\lambda_{32}&=&p_2+p_5+p_6-p_1-p_3-p_4 \\
\end{eqnarray*}

The three correlation matrices $R_1$, $R_2$, $R_3$ are given by:\\
$$
R_1= \left(
\begin{array}{cc}
1& \frac{(p_1+p_2+p_4+p_6-p_3-p_5)-\lambda_{12}\lambda_{13}}{\sqrt{(1-\lambda_{12}^2)(1-\lambda_{13}^2)}} \\
\frac{(p_1+p_2+p_4+p_6-p_3-p_5)-\lambda_{12}\lambda_{13}}{\sqrt{(1-\lambda_{12}^2)(1-\lambda_{13}^2)}} &1\\

\end{array}
\right)
$$
The (1,2) entry of $R_1$ is $R_{23}^1$.
$$
R_2= \left(
\begin{array}{cc}
1& \frac{(p_2+p_3+p_4+p_5-p_1-p_6)-\lambda_{21}\lambda_{23}}{\sqrt{(1-\lambda_{21}^2)(1-\lambda_{23}^2)}} \\
\frac{(p_2+p_3+p_4+p_5-p_1-p_6)-\lambda_{21}\lambda_{23}}{\sqrt{(1-\lambda_{21}^2)(1-\lambda_{23}^2)}} &1 \\

\end{array}
\right)
$$
The (1,2) entry of $R_2$ is $R_{13}^2$.
$$
R_3= \left(
\begin{array}{cc}
1& \frac{(p_1+p_3+p_5+p_6-p_2-p_4)-\lambda_{31}\lambda_{32}}{\sqrt{(1-\lambda_{31}^2)(1-\lambda_{32}^2)}} \\
\frac{(p_1+p_3+p_5+p_6-p_2-p_4)-\lambda_{31}\lambda_{32}}{\sqrt{(1-\lambda_{31}^2)(1-\lambda_{32}^2)}} &1\\

\end{array}
\right)
$$
The (1,2) entry of $R_3$ is $R_{12}^3$.
The probability of a Condorcet winner is
\begin{equation}
L(\delta_{12},\delta_{13},R_{23}^1)+L(\delta_{21},\delta_{23},R_{13}^2)+L(\delta_{31},\delta_{32},R_{12}^3).
\label{eqn10}
\end{equation}
By considering the 27 possibilities with which $\lambda_{12},\lambda_{13}$ and 
$\lambda_{23}$ can take zero, positive and negative values, we can simplify
the expression further. Table \ref{table1} shows the 27 possible cases and
the corresponding limiting value of the probability of a Condorcet winner. A
set of different values is obtained for varying values of $p_i$'s. It is also
to be noted that for a few alternatives the values of 0 and 1 are obtained for
the limiting probability of a Condorcet winner. See alternatives 7 and 17 in Table \ref{table1}. Table \ref{table1} provides the reference table to calculate 
the exact limiting probability for any choice of $p_i,i=1,\cdots,6$. One of the
27 possibilities will apply in a given situation and one can read off
or calculate the value. It is seen from Table \ref{table1} that only in two cases
out of 27, the limiting probability is definitely 0. In 12 cases (more than 40\%)
the probability is 1 and in almost all cases the probability will be at least 50\%. These
observations are consistent with often stated conclusion that for small number
of candidates and large number of votes, the Condercet scheme will produce a
winner.
\begin{table}[t]
\begin{tabular}{|l|l|l|l|l|}
\hline	\hline
Number&$\lambda_{12}$&$\lambda_{13}$&$\lambda_{23}$&Probability of a Winner\\
\hline\hline
1&0&0&0&$\frac{1}{4}+\frac{1}{2\pi}\arcsin(R_{23}^1)+\frac{1}{4}+\frac{1}{2\pi}\arcsin(R_{13}^2)+\frac{1}{4}+\frac{1}{2\pi}\arcsin(R_{12}^3)$\\\hline
2&0&0&+&$\frac{1}{2}+\frac{1}{4} +\frac{1}{2\pi} \arcsin(R_{23}^1)$\\\hline
3&0&0&-&$\frac{1}{2}+\frac{1}{4} +\frac{1}{2\pi} \arcsin(R_{23}^1)$\\\hline
4&0&+&0&$\frac{1}{2}+\frac{1}{4} +\frac{1}{2\pi} \arcsin(R_{13}^2)$\\\hline
5&0&-&0&$\frac{1}{2}+\frac{1}{4} +\frac{1}{2\pi} \arcsin(R_{13}^2)$\\\hline
6&0&+&-&$\frac{1}{2}$\\\hline
7&0&+&+&1\\\hline
8&0&-&+&$\frac{1}{2}$\\\hline
9&0&-&-&1\\\hline
10&+&0&+&$\frac{1}{2}$\\\hline
11&+&+&0&1\\\hline
12&+&0&-&1\\\hline
13&+&-&0&$\frac{1}{2}$\\\hline
14&+&+&+&1\\\hline
15&+&-&-&1\\\hline
16&+&+&-&1\\\hline
17&+&-&+&0\\\hline
18&+&0&0&$\frac{1}{2}+\frac{1}{4} +\frac{1}{2\pi} \arcsin(R_{12}^3)$\\\hline
19&-&0&0&$\frac{1}{2}+\frac{1}{4} +\frac{1}{2\pi} \arcsin(R_{12}^3)$\\\hline
20&-&0&+&1 \\\hline
21&-&0&-&$\frac{1}{2}$\\\hline
22&-&+&0&$\frac{1}{2}$\\\hline
23&-&+&+&1\\\hline
24&-&+&-&0\\\hline
25&-&-&0&1\\\hline
26&-&-&+&1\\\hline
27&-&-&-&1\\
\hline\hline
\end{tabular}
\caption{Limiting Condorcet Winner Probabilities for $m=3$ candidates}
\label{table1}
\end{table}

It can be verified that case 7, for example, is reached when $p_1=p_2=p_5=
\frac{1}{6}$; $p_3=\frac{2}{5}$ and $p_6=\frac{1}{10}$. Case 17 is reached when
$p_1=\frac{3}{22}$,$p_2=\frac{3}{22}$,$p_3=\frac{9}{44}$,$p_4=\frac{13}{44}$,
$p_5=\frac{5}{22}$ and $p_6=0$. For another set of examples for $m=3$, see May \cite{may71}, Garman et al \cite{gar68} and Meyer et al \cite{mey98}. 
%It is straightforward to find probabilities for all $m>3$ for which a Condorcet winner probability is is either 0 or 1.
 
For example, the trivial set of probabilities $p_i=1,p_j=0, j \ne i$ will produce
a Condorcet winner for every $n$ and hence the limit is also 1 
as $n \rightarrow \infty$. Choose the probabilities $p_i$ such that 
$C_1C_2C_3C_4\cdots C_m$, $C_2C_3C_1C_4\cdots C_m$, $C_3C_1C_2C_4\cdots C_m$ 
have probabilities $\frac{1}{3}$ each and the rest of the rankings have 
probabilities 0. It is easily verified that $\lambda_{12}=\frac{1}{3}$, 
$\lambda_{13}=\frac{-1}{3}$, $\lambda_{23}=\frac{1}{3}$, 
$\lambda_{1j}=\lambda_{2j}=\lambda_{3j}=1$ for all $j \ge 4$. 
This implies that $\delta_{13}=\infty$, 
$\delta_{12}=\delta_{23}=-\infty, \delta_{1j}=\delta_{2j}=\delta_{3j}=-\infty$, for $j \ge 4$. From equation (\ref{eqnl}), it follows immediately that in every $L$ term, one of the arguments is a $+\infty$. This implies that the
limiting Condorcet winner probability is 0. Thus for every $m$, we can choose $p_i$ so that the limiting Condorcet winner is 0.

Gillett\cite{gil80} conjectured that the probability of a Condorcet winner $P(3,n;p)$ 
for $m=3$ candidates and any $n$, the number of voters, is minimized for
$p_1=p_4=p_5=\frac{1}{3},p_2=p_3=p_6=0$ or $p_1=p_4=p_5=0,p_2=p_3=p_6=
\frac{1}{3}$. Buckley \cite{buc75} proved this conjecture for the special case $n=3$. 
We will show the validity of the conjecture for general $n$.

The probability of candidate $C_1$ being the majority winner is
$$
P\left\{N_1+N_2+N_3>\frac{n}{2} \hbox{~and~} N_1+N_2+N_5>\frac{n}{2}\right\}
$$
$$
=P\left\{N_1+N_2 > \frac{n}{2}-min(N_3,N_5)\right\}
$$
$$
\ge P\left\{N_1+N_2 >\frac{n}{2}\right\},
$$
since $N_3$ and $N_4$ are non-negative random variables. The above expression 
is greater than or equal to 

$P\left\{N_1>\frac{n}{2}\right\}$, again for the same reason with regard to $N_2$.
$N_1$ is distributed as a Binomial variable with $n$ and $p_1$.
$$
P\left(N_1>\frac{n}{2}\right)=1-B(k;n,p_1)
$$
where $B(k;n,p_1)$ is the Binomial cumulative distribution with parameters
$n$ and $p_1$. $k=\frac{n-1}{2}$ if $n$ is odd and $k=\frac{n}{2}$ if $n$ is even.
From the relation of $B(k;n,p_1)$ with Incomplete Beta function,  see 
Feller(1957) page 127, \cite{ful57},
\begin{equation}
P(N_1>\frac{n}{2})=n {{n-1} \choose k } \int_0^{p_1}{t^k{(1-t)}^{n-k-1}dt}
\label{eqnEE}
\end{equation}

Continuing similar arguments for $C_2$ and $C_3$ being majority winners, 
we have
\begin{eqnarray*}
P(3,n;{\bf p})\ge n {{n-1} \choose k} \left[\int_0^{p_1} {t^k{(1-t)}^{n-k-1}dt}+\int_0^{p_4} {t^k{(1-t)}^{n-k-1}dt}+\int_0^{p_5} {t^k{(1-t)}^{n-k-1}dt}\right]
\end{eqnarray*}

Consider minimizing the following function
\begin{eqnarray}
u(p_1,p_4,p_5)&=& \int_0^{p_1} {t^k{(1-t)}^{n-k-1}dt}+\int_0^{p_4} {t^k{(1-t)}^{n-k-1}dt}\nonumber\\
&+&\int_0^{p_5} {t^k{(1-t)}^{n-k-1}dt} \nonumber\\
\hbox{subject to }p_1+p_4+p_5&=&c  \nonumber \\
p_1 \ge 0 , p_4 \ge 0, p_5 \ge 0 \hbox{ where } c&=&1-(p_2+p_3+p_6)
\label{eqnF}
\end{eqnarray}

We proceed to show that an optimal solution to the problem defined by (\ref{eqnF}) is given
by $p_1=p_4=p_5=\frac{c}{3}$.

$u(p_1,p_4,p_5)$ is a separable function in $p_1$, $p_4$ and $p_5$. It can
be verified that 
\begin{equation}
\frac{\partial^2u(p_1,p_4,p_5)}{\partial {p_i}^2}= {p_i}^{k-1}(1-{p_i}^{n-k-2}) \left(k-(n-1)p_i)\right), i=1,4,5
\label{eqnX}
\end{equation}
For $p_i \le \frac{k}{n-1}, i=1,4,5, \frac{\partial^2u}{\partial {p_i}^2} \ge 0 \Longrightarrow u(p_1,p_4,p_5) $ is convex in $p_i$ in the region of $p_i$'s.

In order that $p_i$'s satisfy the constraints of equation(\ref{eqnF}), we
cannot have two of the $p_i$'s are $\ge \frac{k}{n-1}$ with at least one
of them  $> \frac{k}{n-1}$, since their sum will be greater than 
$\frac{2k}{n-1} \ge 1 $ for both $n$ even and $n$ odd. If for any feasible
solution, a single $p_i > \frac{k}{n-1}$, which implies that the other two
$< \frac{k}{n-1}$, the Hessian matrix will have its determinant a 
negative value because of equation (\ref{eqnX}). Note that the Hessian is a 
diagonal matrix. Thus such a solution 
cannot be a local minimum. Therefore we can confine our solution to the region 
in which all the $p_i$'s  $\le \frac{k}{n-1}$. It is easy to verify that 
$p_1=p_4=p_5=\frac{c}{3}$ satisfies the Kuhn-Tucker condition for the 
minimum. Note also that $\frac{c}{3} < \frac{k}{(n-1)}$, since $\frac{3k}{(n-1)}>1$. Since $u(p_1,p_4,p_5)$
is convex in this region, the solution is a global minimum of
$u(p_1,p_4,p_5)$.

Fix now $p_1=p_4=p_5=\frac{c}{3}$ and repeat the whole argument with
$p_2,p_3,p_6$ as variables and consider the minimization problem of minimizing
the following function.
\begin{eqnarray}
w(p_2,p_3,p_6)&=& \int_0^{p_2} {t^k{(1-t)}^{n-k-1}dt}+\int_0^{p_3} {t^k{(1-t)}^{n-k-1}dt}\nonumber\\
&+&\int_0^{p_6} {t^k{(1-t)}^{n-k-1}dt} \nonumber\\
\hbox{subject to }p_2+p_3+p_6&=&1-c  \nonumber \\
p_2 \ge 0 , p_3 \ge 0, p_6 \ge 0 \hbox{ where } c&=&(p_1+p_4+p_5)
\label{eqnG}
\end{eqnarray}

As before, the minimum solution of equation (\ref{eqnG}) is
$p_2=p_3=p_6=\frac{1-c}{3}$ for $p_1=p_4=p_5=\frac{c}{3}$. From equation (\ref{eqnG}), we see that the least value of $w(p_3,p_3,p_6)$ is when $c=1$.
Thus the minimum solution is $p_1=p_4=p_5=\frac{1}{3};~p_2=p_3=p_6=0$.
Note that our arguments show that $p_2=p_3=p_6=\frac{1}{3};~p_1=p_4=p_5=0 $
is also optimal.
Denote by ${\bf{\tilde{p}^*}}$ the solution ${p_1}^*={p_4}^*={p_5}^*=\frac{1}{3}$;
${p_2}^*={p_3}^*={p_6}^*=0$. Let $u^*({p_1}^*,{p_4}^*,{p_5}^*) $ be the value
of the objective function of (\ref{eqnF}). We have shown that for any choice of
$\tilde{p}$
\begin{eqnarray*}
P(3,n;\tilde{p}) &\ge& n {{n-1}\choose k} u^*({p_1}^*,{p_4}^*,{p_5}^*)\\
&=& 3(1-B(k;n,\frac{1}{3})), \hbox{ from (\ref{eqnE}) }
\end{eqnarray*}

We will now show that $P(3,n;{\bf{\tilde{p}^*}})=3(1-B(k;n,\frac{1}{3}))$
which concludes that ${\bf{\tilde{p^*}}}$ yields the minimum value of 
$P(3,n;\tilde{p})$. Since ${p_2}^*={p_3}^*={p_6}^*=0$, we have
\begin{eqnarray*}
P(3,n;{\bf{\tilde{p}^*}})&=&P(N_1 >\frac{n}{2},N_1+N_5 >\frac{n}{2})\\
&+&P(N_4 >\frac{n}{2},N_4+N_1 >\frac{n}{2})\\
&+&P(N_5 >\frac{n}{2},N_5+N_4 >\frac{n}{2})\\
&=& P(N_1>\frac{n}{2}) + P(N_4>\frac{n}{2})+P(N_5>\frac{n}{2})\\
&=&3(1-B(k;n,\frac{1}{3}),
\end{eqnarray*}
since $N_1,N_4,N_5$ are identically Binomial distributed with parameters
$n$ and $\frac{1}{3}$.

This establishes the conjecture. Note that the minimum value for the
probability of a Condorcet winner is given by 
$$
3\left(1-B(k;n,\frac{1}{3})\right)
$$
or an equivalent integral representation
$$
3n{{n-1}\choose k} \left[\int_0^{\frac{1}{3}} {t^k{(1-t)}^{n-k-1}dt}\right]
$$
The former is easier to evaluate for small values of $n$ and employ
a normal approximation to the Binomial for large $n$. For example, when $n=3$,
the minimum probability is
$$
3\left({3 \choose 2} { \left( \frac{1}{3}\right)} ^2 \left(\frac{2}{3}\right)+{3 \choose 3} {\left(\frac{1}{3}\right)}^3 \right)=\frac{7}{9}
$$
which is about 77\%. For $n=4$, the minimum probability of a Condorcet
winner is
$$
3\left({4 \choose 3}{ \left(\frac{1}{3}\right)}^3\left(\frac{2}{3}\right)+{4 \choose 4} {\left(\frac{1}{3}\right)}^4 \right)=\frac{1}{3}
$$
which drops down to 33\%. 

The foregoing arguments can be generalized to $m>3$ candidates. Consider the
probability distribution $\bf{\tilde{p}}^*$ on $m!$ preference rank orders given by assigning
probability of $\frac{1}{m}$ to each of the cyclic rank orders 
\begin{eqnarray}
\hbox{ rank order } 1,2,\cdots,(m-1),m &~& \hbox {probability } p_1 \nonumber\\
\hbox{  rank order } 2,3,\cdots,m,1 &~&\hbox{ probability } p_2 \nonumber \\
\hbox{ rank order } 3,4,\cdots,1,2 &~& \hbox { probability } p_3 \nonumber \\
\cdots&~& \cdots \nonumber \\
\hbox{ rank order } m,1,\cdots, (m-2),(m-1)&~& \hbox{ probability } p_m
\label{eqnAAA}
\end{eqnarray}
We will establish that 
this solution yields the minimum for the probability of a Condorcet winner 
$P(m,n;\tilde{p})$ for a general probability distribution $\tilde{p}$. We use
the same approach that was used for $m=3$ case. It is shown first along
similar lines that
\begin{eqnarray}
P(m,n;\tilde{p}) &\ge& P(N_1>\frac{n}{2}) + P(N_2>\frac{n}{2})+\cdots+P(N_m>\frac{n}{2})\\
&=& n {{n-1} \choose k} \left[ \displaystyle\sum_{i=1}^m{\int_0^{p_i}{t^k{(1-t)}^{n-k-1} dt}} \right]
\end{eqnarray}
From an argument similar to the case $m=3$, it is shown that
\begin{equation}
P(m,n;\tilde{p}) \ge m (1-B(k;n,\frac{1}{m})).
\label{eqnBB}
\end{equation}

Next it follows that from (\ref{eqnAAA}) and the definition of a Condorcet winner,
$P(m,n;{\bf{\tilde{p}}^*})$ for 
${\bf{\tilde{p}^*}}$ given by 
${\bf{p_i}^*}=\frac{1}{m},$ for $i=1,\cdots,m$ 
and ${\bf{p_i}^*}=0$ for
$i \ge (m+1)$ the probability of $C_1$ winning is
\begin{eqnarray*}
P(N_1+N_3+\cdots+N_m >\frac{n}{2},N_1+N_4+\cdots+N_m >\frac{n}{2},\\
\cdots, N_1> \frac{n}{2})\\
&=&P(N_1>\frac{n}{2})
\end{eqnarray*}
Similarly it is verified that the probability of $C_2$ winning is
\begin{eqnarray*}
P(N_1+N_2+N_4+\cdots+N_m >\frac{n}{2},N_1+N_2+N_5+\cdots+N_m >\frac{n}{2},\\
\cdots, N_2> \frac{n}{2})\\
&=&P(N_2>\frac{n}{2})
\end{eqnarray*}
and so on. Hence
\begin{eqnarray*}
P(m,n;{\bf{\tilde{p}}^*}) &=& P(N_1>\frac{n}{2}) + P(N_2>\frac{n}{2})+\cdots+P(N_m>\frac{n}{2})\\
&=&m (1-B(k;n,\frac{1}{m}))
\end{eqnarray*}
which equals the bound in (\ref{eqnBB}). Hence ${\bf{\tilde{p}^*}}$ yields the 
minimum for the probability of a Condorcet winner.
Table \ref{table2} gives the minimum probability of a
Condorcet winner for selected values of $m$ and $n$.

\begin{table}[t]
\begin{center}
\begin{tabular}{|l|l|l|l|l|}
\hline	\hline
$n$&$m=3$&$m=4$&$m=5$&$m=10$\\
\hline\hline
3&       0.7778 & 0.6250&   0.5200&    0.2800\\\hline
4&       0.3333&  0.2031&  0.1360&    0.0370\\\hline
5&       0.6296&  0.4141&  0.2896&    0.0856\\\hline
6&       0.3004&  0.1504&  0.0848&    0.0127\\\hline
7&       0.5199&  0.2822&  0.1667&    0.0273\\\hline
8&       0.2638&  0.1092&  0.0520&    0.0043\\\hline
9&       0.4345&  0.1957&  0.0979&    0.0089\\\hline
10&      0.2297&  0.0789&  0.0318&   0.0015\\\hline
19&      0.1943&  0.0356&  0.0079&   0.0000\\\hline
20&      0.1129&  0.0158&  0.0028&   0.0000\\\hline
29&      0.0934&  0.0071&  0.0007&   0.0000\\\hline
30&      0.0564&  0.0033&  0.0003&   0.0000\\\hline
50&      0.0148&  0.0002&  0.0000&   0.0000\\\hline
51&      0.0207&  0.0002&  0.0000&   0.0000\\\hline
100&     0.0006&  0.0000&  0.0000&   0.0000\\\hline
101&     0.0008&  0.0000&  0.0000&   0.0000\\\hline\hline
\end{tabular}
\caption{Minimum Condorcet Winner Probabilities for some values of $m$ and $n$}
\label{table2}
\end{center}
\end{table}

\section{Condorcet Winner: Special Case m=4 candidates}
For this case, we show that it is possible to obtain, for any choice of $p_i,i=1,2,\cdots,24$, 
exact limiting probability of a Condorcet winner. We first calculate the six 
quantities $\lambda_{12},\lambda_{13},\lambda_{14},\lambda_{23},\lambda_{24}$ and
$\lambda_{34}$ from (\ref{eqn49}). Then define
$\delta_{12},\delta_{13},\delta_{14},\delta_{23},\delta_{24}$ and $\delta_{34}$
as stated in Section 2.
Also calculate the correlation matrices $R_i,i=1,2,3,4$. In this case,
there are $3^6=729$ possibilities for the $\lambda_{ij}$'s to be $>0$ or $<0$ or $=0$.
It will be unwieldy to prepare a Table containing all these 729 cases.
We show how the limiting probability can be calculated for any given
situation. The limiting probability of Condorcet winner follows from (\ref{eqnll}).
\begin{eqnarray*}
P(\infty,4)&=& L(\delta_{12},\delta_{13},\delta_{14};R_1)+L(\delta_{21},\delta_{23},\delta_{24};R_2)\\
&~&+L(\delta_{31},\delta_{32},\delta_{34};R_3)+L(\delta_{41},\delta_{42},\delta_{43};R_4)
\end{eqnarray*}
Consider the calculation of $L(\delta_{12},\delta_{13},\delta_{14};R_1)$. 
Calculations for other L-factors are similar. To  obtain this, we need at the
worst to calculate $L(0,0,0;R_1)$ or $L(0,0,R_{ij}^1)$ or $L(0)$.
Recall the definition of $R_{ij}^1 $ as appropriate $i,j$ element of $R_1$.
Note that $L(0)=\Phi(0)=\frac{1}{2}$, where $\Phi$ is the cdf of the standard
normal distribution. Bacon \cite{bac63} gives exact expressions for those $L$ functions as
\begin{eqnarray*}
L(0,0,R_{ij}^1)&=& \frac{1}{4}+\frac{1}{2\pi} arcsin(R_{ij}^1)\\
L(0,0,0,R_1)&=& \frac{1}{8}\left[1+\frac{2}{\pi} \left(arcsin(R_{23}^1)
+arcsin(R_{24}^1)+arcsin(R_{34}^1)\right)\right]
\end{eqnarray*}
The calculation of this $L$ for other possible arguments is by inspection or
it will reduce to one of the above 3 cases with arguments equal to zero.

This example illustrates the comment made immediately after equation (\ref{eqnll})
that the only calculations, if any, involve calculating the positive orthant probabilities of
appropriate multivariate normal distribution. Thus for any one of the 729 possibilities
that is reached in a given case, we can calculate the exact limiting
probability of a Condorcet winner.
\section{Conclusion}
In this paper, we provide a rigorous mathematical treatment for calculating
the limiting probability of a Condorcet winner for any number of candidates
and with arbitrary rank order probabilities, when the voters are independent.
We show further that the limiting probability depends only on the
positive orthant probabilities of an appropriate multinormal distribution
even for arbitrary probabilities of the preference rank order
probabilities.
We provide approximate
limiting probabilities for a Condorcet winner for any number of candidates and
equal rank order probabilities. We provide a compact and complete
Table for the limiting probability of a Condorcet winner with 3 candidates 
and arbitrary rank order probabilities. We present a simple proof of a result
of May \cite{may71} to show that the limiting probability of a Condorcet 
winner tends to zero as the number of candidates tends to infinity. We
also present a scheme for calculating the exact limiting probability  for
$m=4$ candidates and arbitrary probability distributions on the preference
rank orders. We show for the first time that the limiting probability of a
Condorcet winner for any given number of candidates $m$ is monotone decreasing in
$m$ for the equally likely case. This, in turn, settles the conjectures of 
Kelly \cite{kel74} and Buckley and Westen \cite{bucwes79} for the case $n \rightarrow \infty$.
We prove the validity of Gillett's \cite{gil80} conjecture on the minimum value
of the probability of a Condorcet winner for $m=3$ and any $n$. Extending the
result we obtain the minimum value and the minimum solution for any $n$ and $m$.
\bibliography{condorcet}
\end{newpage}
\end{newpage}
\end{document}